\documentclass[12pt,reqno]{amsart}
\usepackage{enumerate, latexsym, amsmath, amsfonts, amssymb, amsthm, graphicx, color}
\usepackage{hyperref}
\hypersetup{hypertex=true,
	colorlinks=true,
	linkcolor=blue,
	anchorcolor=blue,
	citecolor=blue}
\usepackage{amsmath}
\def\pmod #1{\ ({\rm{mod}}\ #1)}
\def\i{\mathrm{i}}
\def\Z{\Bbb Z}

\def\C{\Bbb C}

\def\l{\left}
\def\r{\right}
\def\bg{\bigg}
\def\({\bg(}
\def\){\bg)}
\def\t{\text}
\def\f{\frac}

\def\ls{\le}

\def\ve{\varepsilon}

\def\eq{\equiv}

\theoremstyle{plain}
\newtheorem{theorem}{Theorem}

\newtheorem{lemma}{Lemma}

\theoremstyle{definition}

\theoremstyle{remark}

\newcommand{\legendre}[2]{\genfrac(){}{}{#1}{#2}}
\newcommand{\eps}{\varepsilon}
\makeatletter
\@namedef{subjclassname@2010}{%
	\textup{2010} Mathematics Subject Classification}
\makeatother
\vspace{4mm}

\allowdisplaybreaks[1]
\begin{document}	
	\hbox{Preprint}
\medskip

	\title[Matrices composed of Legendre symbols]
{Symmetric identities and determinants for matrices composed of Legendre symbols}
	
	\author[R. Le Lan and K. Liu]{Roman Le Lan and Keqin Liu*}

    \address {(Roman Le Lan) Lyc\'ee Louis-Le-Grand, 75005, Paris, France} \email{roman.lelan1@gmail.com} 
	
	\address {(Keqin Liu) School of Mathematics and Physics, Xi'an Jiaotong-Liverpool University,
	Suzhou, 215123, People's Republic of China} \email{Keqin.Liu@xjtlu.edu.cn}

	\begin{abstract}
   In this paper, we prove a series of identities and evaluate the determinants for a class of matrices composed of Legendre symbols given by Sun \cite[Conjecture 3.2 (i)]{S24}. Together with \cite{LSW} which proved Conjecture~3.2~(ii) in \cite{S24}, these results provide a more general approach for a complete proof of Conjecture~3.2 in \cite{S24}.
	\end{abstract}
	
	\thanks{2020 {\it Mathematics Subject Classification}.
		Primary 11A15, 11C20; Secondary 15A15.
		\newline\indent {\it Keywords}. Determinant, Legendre symbol, quadratic field.
		\newline \indent
		The work of Keqin Liu was supported by Leadership Talent Program (Science and Education) of SIP (KJQ2024202) and XJTLU Research Development Fund (RDF-24-01-048).
		\newline\indent{*Corresponding author.}
        }
	\maketitle
	\section{Introduction}	
	\setcounter{lemma}{0}
	\setcounter{theorem}{0}
	\setcounter{corollary}{0}
	\setcounter{remark}{0}
	\setcounter{equation}{0}
	\setcounter{conjecture}{0}

Let $\det A$ or $|a_{ij}|_{1\ls i,j\ls n}$ denote the determinant of a matrix $A=[a_{ij}]_{1\ls i,j\ls n}$ over a field. Let $p>2$ be an odd prime and~$(\f{\cdot}p)$ the Legendre symbol. For $p\equiv 1\pmod4$, let $\ve_p>1$ and $h_p$ denote the fundamental unit and class number of the real quadratic field $\mathbb{Q}(\sqrt{p})$, respectively. For convenience, we write
\begin{align}\label{abp}
	 \varepsilon_p^{h_p}=a_p+b_p\sqrt{p}\ \ \t{with}\ a_p,b_p\in\mathbb{Q},
\end{align}
and
\begin{align}\label{ab'p}
	 \varepsilon_p^{(2-(\frac{2}{p}))h_p}=a_p'+b_p'\sqrt{p}\ \ \t{with}\ a_p',b_p'\in\mathbb{Q}.
\end{align}
In 2024, Sun~\cite{S24} made a series of conjectures on the determinants of various matrices represented by linear combinations of Legendre symbols, including the following one.
\begin{theorem}[Conjecture 3.2 in~\cite{S24}]\label{thm:main}
Let~$p$ be an odd prime, and set $v = wx - (y + 1)(z + 1)$.

{\rm (i)} If $p>3$, then
\begin{equation}\label{main1}
 \begin{aligned}&\ \l|x+\l(\f{i+j}p\r)+\l(\f ip\r)y+\l(\f jp\r)z+\l(\f{ij}p\r)w\r|_{0\ls i,j\ls(p-3)/2}
\\=&\ \begin{cases} \l(\f 2p\r)2^{\frac{p-3}{2}}((pb_p-2a_p)x+(a_p-2b_p)v)&\t{if}\ p\eq1\pmod4,
\\2^{\frac{p-3}{2}}(v-2x)&\t{if}\ p\eq3\pmod4.
\end{cases}
\end{aligned}
\end{equation}

{\rm (ii)} We have
\begin{equation}\label{main2}
 \begin{aligned}&\ \l|x+\l(\f{i-j}p\r)+\l(\f ip\r)y+\l(\f jp\r)z+\l(\f{ij}p\r)w\r|_{0\ls i,j\ls(p-3)/2}
\\=&\ \begin{cases} -a_p'x - \l(\f 2p\r)b_p'v&\t{if}\ p\eq1\pmod4,
\\x&\t{if}\ p\eq3\pmod4.
\end{cases}
\end{aligned}
\end{equation}

\end{theorem}
Case~(ii) in Theorem~\ref{thm:main} was proven in \cite{LSW}. In this paper, we prove Case~(i) above. Our approach is based on the evaluation of certain quadratic and bilinear forms involving the $p$-th root of unity $\zeta:=e^{2\pi \i/p}$, where~$e$ is Euler's number and~$\i$ the imaginary unit in~$\mathbb{C}$. 

\subsection*{Notations} 
We will use the following notations. Let $p>3$ be a prime number. Define
\[c_p:=\prod_{k=1}^{(p-1)/2}\left(1-\zeta^{k^2}\right)\]
and
\[N:=\left\{k\mod{p\Z}: \legendre{k}{p}=-1\right\}.\]
Furthermore, define $N_0=N\cup\{0\}$. Note that $c_p$ is determined for all $p>3$ by Sun \cite[Theorem~1.3]{SunFFA}. Let $\tau_p$ denote the quadratic Gauss sum \cite{IK}: 
\begin{equation}\label{GaussSum}
\tau_p=1+2\sum_{k=1}^{(p-1)/2}\zeta^{k^2}=\begin{cases} \sqrt{p}&\t{if}\ p\eq1\pmod4,\\\i \sqrt{p}&\t{if}\ p\eq3\pmod4.
\end{cases}
\end{equation}
For any integer $m$ not divisible by $p$, by the Galois theory, we have
\begin{equation}\label{eqn:gauss-galois}
\legendre{m}{p}=\frac{1}{\tau_p}\left(1+2\sum_{k=1}^{(p-1)/2}\zeta^{mk^2}\right).
\end{equation}
\section{Proof of Theorem~\ref{thm:main}~(i)}
	\setcounter{lemma}{0}
	\setcounter{theorem}{0}
	\setcounter{corollary}{0}
	\setcounter{remark}{0}
	\setcounter{equation}{0}
	\setcounter{conjecture}{0}

In this section, we prove the identity \eqref{main1} of Theorem ~\ref{thm:main}.
\subsection{Reduction}

Let $x,y,z,w$ be complex numbers and $p>3$ a prime. Define
\[A(x,y,z,w):=\left\lvert x+\legendre{i+j}{p}+\legendre ipy+\legendre jpz+\legendre{ij}pw\right\rvert_{0\le i,j\le (p-3)/2}\]
and
\[H(x):=\left\lvert x+\legendre{i+j}{p}\right\lvert_{0\le i,j\le (p-3/2}.\]
By \cite[Lemma 2.1]{RJ},
\begin{equation*}
    H(x)=H(0)+x\left\lvert\legendre{i+j}p-\legendre ip-\legendre jp\right\rvert_{1\le i,j\le (p-3)/2}.
\end{equation*}
Then, by \cite[(1.13)]{S24},
\begin{align*}
A(x,y,z,0)&=(y+1)(z+1)H(0)+x\left\lvert\legendre{i+j}p-\legendre ip-\legendre jp\right\rvert_{1\le i,j\le (p-3)/2}\\&=(y+1)(z+1)H(0)+H(x)-H(0)\\&=(yz+y+z)H(0)+H(x).
\end{align*}
Also, by \cite[(1.11)]{S24},
\[A(0,y,z,w)=A(0,y,z,0).\]
Therefore, we find that
\[A(1,0,0,0)=H(1),\quad A(0,1,0,0)=A(0,0,1,0)=2H(0),\]
\[\text{and}\quad A(0,0,0,1)=H(0).\]
Finally, by \cite[Theorem 3.1]{ChenSun},
\begin{equation}\label{reduction}
A(x,y,z,w)=xH(1)+(1-x+y+z+yz-wx)H(0).
\end{equation}

\begin{lemma}[Generalized Matrix Determinant Lemma]\label{lma:matrix}
Let $u,v$ be two $m\times1$ column vectors, $a,b\in\C$ and $A$ an $m\times m$ invertible matrix. Then,
\[\det(A+auu^T+bvv^T)=\det(A)\left((1+a\alpha)(1+b\beta)-ab\gamma\delta\right),\]
where
\[\alpha=u^TA^{-1}u,\quad \beta=v^TA^{-1}v,\quad \gamma=u^TA^{-1}v,\quad \delta=v^TA^{-1}u.\]
\end{lemma}
\begin{proof}
We prove the lemma by using the following generalization of the Matrix Determinant Lemma \cite{mdl,wiki}:
\[\det(A+UWU^T)=\det(I+WU^TA^{-1}U)\det(A),\]
where $U$ is any $m\times r$ matrix, $W$ any $r\times r$ matrix, and $I$ the identity matrix. Taking $U=[u,v]$ and
\[W=\begin{pmatrix}
       a&0\\0&b
\end{pmatrix},\]
we arrive at the desired identity immediately.
\end{proof}

Let $p=2n+1$ and
\[H(x):=\left\lvert x+\legendre{i+j}{p}\right\lvert_{0\le i,j\le n-1}.\]
By the identity \eqref{eqn:gauss-galois}, we have 
\[\tau_p\left(x+\legendre{i+j}p\right)=\tau_px+1+2\sum_{k=0}^{n-1}\zeta^{i(k+1)^2}\zeta^{j(k+1)^2},\]
except for $i=j=0$. Define 
\[M:=\left(\zeta^{i(j+1)^2}\right)_{0\le i,j\le n-1}.\]
Set $u=[1,1,\ldots,1]^T$ be the $n\times 1$ column vector with all $1$s and $v=[1,0,\ldots,0]^T$ the $n\times 1$ column vector with the first entry equal to $1$ and the rest all equal to $0$. Then,
\begin{align}\nonumber
\tau_p^nH(x)=&\det\left(2MM^T+(x\tau_p+1)uu^T-(2n+1)vv^T\right)\\\nonumber
=&\det\left(2MM^T+(x\tau_p+1)uu^T-pvv^T\right).
\end{align}
Substituting $A=2MM^T$ in Lemma~\ref{lma:matrix}, we have
\begin{align}\label{eqn:reduce}
\tau_p^nH(x)=&2^{n}\det(M)^2\\\nonumber&\cdot\left[\left(1+(\tau_px+1)\frac12\alpha_1\right)\left(1-\frac p2\beta_1\right)+(\tau_px+1)\frac p4\gamma_1\delta_1\right],
\end{align}
where
\begin{align*}
\alpha_1&=(M^{-1}u)^TM^{-1}u, \quad \beta_1=(M^{-1}v)^TM^{-1}v,\\
\gamma_1&=(M^{-1}u)^TM^{-1}v, \quad \delta_1=(M^{-1}v)^TM^{-1}u.
\end{align*}
Note that the invertibility of $M$ (and thus~$A$) will be established below. 

\subsection{Evaluation of $\det(M)^2$}
By the Vandermonde determinant formula, we have \cite{ST}:
\begin{align*}
\det(M)=\prod_{1\le i<j\le n}(\zeta^{j^2}-\zeta^{i^2}).
\end{align*}
Then, by \cite[Theorem~1.3]{SunFFA}, we have 
\begin{equation}\label{detM}
\det(M)^2=\begin{cases} (-1)^{\frac{p-1}4}p^{\frac{p-3}4}\eps_p^{h_p}=\legendre{2}{p}p^{\frac{p-3}4}\eps_p^{h_p}\neq0&\t{if}\ p\eq1\pmod4,\\
(-p)^{\frac{p-3}4}\neq0&\t{if}\ p\eq3\pmod4.
\end{cases}
\end{equation}

\subsection{Evaluation of $\alpha_1$}
The evaluations of $\alpha_1$, $\beta_1$, $\gamma_1$ and $\delta_1$ rely on the following lemma.
\begin{lemma}[Partial Fraction Decomposition]\label{lma:pfd}
    Let $P$ and $Q$ be two polynomials such that 
    \[Q=\prod_{k=1}^m(X-\lambda_k)\]
    where the $\lambda_k$ are pairwise distinct for all $1\le k\le m$.
    Then
    \[\frac PQ= S+\sum_{k=1}^m\frac{P(\lambda_k)}{Q'(\lambda_k)(X-\lambda_k)},\]
    where $S$ is the quotient in the Euclidean division of $P$ by $Q$, and $'$ denotes the derivative with respect to $X$.
\end{lemma}
\begin{proof}
    By Euclidean division, there exists a polynomial $R$ of degree lower than $\deg Q$ and a polynomial $S$ such that
    \[P=SQ+R.\]
    Then, by Lagrange's interpolation \cite{V3}, we can factorize $R/Q$ into a linear combination of distinct terms $\{1/(X-\lambda_k)\}_{1\le k\le m}$:
    \[\frac{P}{Q}=S+\frac RQ=S+\sum_{k=1}^m \frac{R(\lambda_k)}{Q'(\lambda_k)(X-\lambda_k)}.\]
    The lemma follows immediately by noting that $P(\lambda_k)=R(\lambda_k)$ for all $1\le k\le m$.
\end{proof}

Let us represent the quadratic form $(M^{-1}u)^TM^{-1}u$ as the sum of squares:
\[\alpha_1=(M^{-1}u)^TM^{-1}u=\sum_{k=1}^{n}a_k^2,\]
where $M^{-1}u=[a_1,\ldots,a_n]^T$. By matrix multiplication, for all $0\le i\le n-1$,
\[\sum_{k=1}^{n}a_k\zeta^{ik^2}=1.\]
Consider an arbitrary polynomial with degree at most $n-1$: $P=\sum_{i=0}^{n-1}d_iX^i$, where $\left(d_0,\ldots,d_{n-1}\right)$ are complex numbers. Then,
\[\sum_{k=1}^na_kP\left(\zeta^{k^2}\right)=P(1).\]
Replacing $P$ by the $k$-th Lagrange polynomial $L_k(X)$ for Lagrange interpolation at nodes $\left(\zeta,\ldots,\zeta^{n^2}\right)$ and noting that $L_k\left(\zeta^{j^2}\right)=\delta_{kj}$ with $\delta$ the Kronecker symbol \cite{V3}, we have
\[a_k=L_k(1)=\prod_{\substack{1\le j\le n\\j\ne k}}\frac{1-\zeta^{j^2}}{\zeta^{k^2}-\zeta^{j^2}}.\]
Hence, using the identity \cite[(2.8)]{WW},
\begin{align*}
\alpha_1&=\sum_{k=1}^n\prod_{\substack{1\le j\le n\\j\ne k}}\frac{\left(1-\zeta^{j^2}\right)^2}{\left(\zeta^{k^2}-\zeta^{j^2}\right)^2}
\\&=-\frac{c_p^2}p\sum_{k=1}^n\frac{\zeta^{k^2}\prod_{i\in N}\left(\zeta^{k^2}-\zeta^{i}\right)}{1-\zeta^{k^2}}\prod_{\substack{1\le j\le n\\j\ne k}}\frac{1}{\zeta^{k^2}-\zeta^{j^2}}.
\end{align*}
Let
\[P_1=X\prod_{i\in N}\left(X-\zeta^{i}\right)\quad\text{and}\quad Q_1=\prod_{k=1}^n\left(X-\zeta^{k^2}\right).\]
By coefficient comparison \cite{Liu25LMA}, one can easily obtain the Euclidean division of $P_1$ by $Q_1$:
\[P_1=(X+\tau_p)Q_1+R,\]
with $R$ a polynomial of degree lower than $\deg Q_1$. Using Lemma~\ref{lma:pfd} at $X=1$, we have
\[\alpha_1=-\frac{c_p^2}p\left(\frac{P_1(1)}{Q_1(1)}-1-\tau_p\right)\]
Since 
\[P_1Q_1=X\frac{X^p-1}{X-1},\]
as $X$ tends to $1$, one can get, after some simplifications,
\[\alpha_1=\frac{c_p^2}p\left(1+\tau_p\right)-1.\]

\subsection{Evaluation of $\beta_1$}
We use a method similar to the one used for the evaluation of $\alpha_1$. Let us represent the quadratic form $(M^{-1}v)^TM^{-1}v$ as the sum of squares:
\[\beta_1=(M^{-1}v)^TM^{-1}v=\sum_{k=1}^{n}b_k^2,\]
where $M^{-1}v=[b_1,\ldots,b_n]^T$. By matrix multiplication, for all $0\le i\le n-1$,
\[\sum_{k=1}^{n}b_k\zeta^{ik^2}=\begin{cases}
    1&\text{if $i=0$}\\ 0&\text{else.}
\end{cases}.\]
Then, by evaluating the $k$-th Lagrange polynomial $L_k(X)$ at $X=0$ with $\deg L_k\le n-1$, one can get
\[b_k=\prod_{\substack{j=1\\j\ne k}}^n\frac{-\zeta^{j^2}}{\zeta^{k^2}-\zeta^{j^2}}.\]
Hence, since $\prod_{k=1}^n\zeta^{k^2}=1$, using the identity \cite[(2.8)]{WW} again, we have
\begin{align*}
\beta_1&=\sum_{k=1}^n\prod_{\substack{j=1\\j\ne k}}^n\frac{\zeta^{2j^2}}{\left(\zeta^{k^2}-\zeta^{j^2}\right)^2}\\&=\sum_{k=1}^n\zeta^{-2k^2}\prod_{\substack{j=1\\j\ne k}}^n\frac{1}{\left(\zeta^{k^2}-\zeta^{j^2}\right)^2}\\
&=\frac1p\sum_{k=1}^n\zeta^{-k^2}\prod_{i\in N_0}\left(\zeta^{k^2}-\zeta^i\right)\prod_{\substack{j=1\\j\ne k}}^n\frac{1}{\zeta^{k^2}-\zeta^{j^2}}.
\end{align*}
Let
\[P_2=\prod_{i\in N_0}\left(X-\zeta^{i}\right)\quad\text{and}\quad Q_2=\prod_{k=1}^n\left(X-\zeta^{k^2}\right).\]
It can be easily established that the quotient of the Euclidean division of $P_2$ by $Q_2$ is $X+\tau_p-1$. 
Then, using Lemma~\ref{lma:pfd} at $X=0$, after some simplifications, we get
\[\beta_1=\frac{\tau_p}{p}.\]

\subsection{Evaluation of $\gamma_1$ and $\delta_1$}
Once again, we are proceeding in the same way above. Clearly,
\[\gamma_1=\delta_1=\sum_{k=1}^na_kb_k.\]
Then, using the identity \cite[(2.8)]{WW} again,
\begin{align*}
\gamma_1=\delta_1&=\sum_{k=1}^na_{k}b_{k}=\sum_{k=1}^n\prod_{\substack{j=1\\j\ne k}}^n\frac{-\zeta^{j^2}\left(1-\zeta^{j^2}\right)}{\left(\zeta^{k^2}-\zeta^{j^2}\right)^2}\\
&=\frac{(-1)^{n-1}}p\sum_{k=1}^n\prod_{i\in N_0}\left(\zeta^{k^2}-\zeta^{i}\right)\prod_{\substack{j=1\\j\ne k}}^n\frac{1-\zeta^{j^2}}{\zeta^{k^2}-\zeta^{j^2}}\\
&=\frac{(-1)^{n-1}c_p}p\sum_{k=1}^n\frac{\prod_{i\in N_0}\left(\zeta^{k^2}-\zeta^{i}\right)}{1-\zeta^{k^2}}\prod_{\substack{j=1\\j\ne k}}^n\frac{1}{\zeta^{k^2}-\zeta^{j^2}}
\end{align*}
Using the same partial fraction decomposition as for evaluating $\beta_1$, but evaluating the polynomials at $X=1$, we get
\[\gamma_1=\delta_1=\frac{(-1)^{n-1}c_p}p\left(\frac{P_2(1)}{Q_2(1)}-\tau_p\right)=\frac{(-1)^nc_p}p\tau_p.\]

\subsection{Final evaluation}
Substituting the values of $\alpha_1$, $\beta_1$, $\gamma_1$, and $\delta_1$ into \eqref{eqn:reduce}, after some simplifications, we obtain
\begin{equation}\label{eqn:Hx}
    H(x)=2^{n}\det(M)^2\cdot\frac{2c_p^2\tau_px+2c_p^2+c_p^2\tau_p^2x+c_p^2\tau_p+2p-2p\tau_px-p\tau_p+p\tau_p^2x}{\tau_p^{n}4p}.
\end{equation}
By \cite[Theorem~1.3]{SunFFA}, and noting that the norm of $\varepsilon_p$ with respect to $\mathbb{Q}(\sqrt{p})/\mathbb{Q}$ is $-1$ and the class number $h_p$ is odd when $p\eq1\pmod4$  (cf. \t{\cite[pp. 185-187]{Cohn}} or \t{\cite[Pg.~237]{chapman}}), we have
\begin{equation}\label{cp}
c_p=\begin{cases} \sqrt{p}\ve_p^{-h_p} = -a_p\sqrt{p}+b_pp&\t{if}\ p\eq1\pmod4,\\
(-1)^{(h_{-p}+1)/2}\i\sqrt{p}&\t{if}\ p\eq3\pmod4,
\end{cases}
\end{equation}
where $h_{-p}$ is the class number of the imaginary quadratic field $\mathbb{Q}(\sqrt{-p})$ and is odd when $p\eq3\pmod4$.

By~\eqref{GaussSum}, \eqref{detM} and~\eqref{cp}, after some simplifications, we arrive at
\begin{equation}\label{HxCF}
    H(x)=\begin{cases} \legendre{2}{p}2^{\frac{p-3}{2}}((pb_p-2a_p)x+(2b_p-a_p))&\t{if}\ p\eq1\pmod4,\\
-2^{\frac{p-3}{2}}(1+2x)&\t{if}\ p\eq3\pmod4,
\end{cases}
\end{equation}
Substituting~\eqref{HxCF} into~\eqref{reduction}, we arrive at~\eqref{main1} immediately.

\section{Concluding Remarks}
Another example of applying our approach by reducing the matrix determinants to evaluations of constants determined by symmetric polynomials is given by Chapman's original matrix determinant \cite{chapman}:
\[F(x):=\left|x+\legendre{i+j}{p}\right|_{1\le i,j\le n}.\]
Define
\[M_1:=\left(\zeta^{ij^2}\right)_{1\le i,j\le n}.\]
Then
\begin{align*}
    \tau_p^{n}F(x)=2^{n-1}\det(M_1)^2\left(2+(\tau_px+1)(M_1^{-1}u)^TM_1^{-1}u\right).
\end{align*}
Using the method above for evaluating the quadratic form $(M_1^{-1}u)^TM_1^{-1}u$, we have an alternative proof of \cite[$\det C_p(x)$ in Corollary~3]{chapman}.

\medskip
\noindent{\bf Statements and Declarations}. There are no competing interests. This original paper contains no data, and  it has not been submitted elsewhere.


\end{document}